# On the Uniqueness of Participation Factors in Nonlinear Dynamical Systems

Tianwei Xia and Kai Sun*
Department of Electrical Engineering and Computer Science,
University of Tennessee, Knoxville, USA
kaisun@utk.edu

**Abstract:** In the modal analysis and control of nonlinear dynamical systems, the participation factors of state variables with respect to a critical or selected mode serve as a pivotal tool for simplifying stability studies by focusing on a subset of highly influential state variables. For linear systems, the participation factors of state variables regarding a mode are uniquely determined by the mode's composition and shape, defined by the system's left and right eigenvectors, respectively. However, the uniqueness of other types of participation factors necessitates further investigation. This paper establishes a sufficient condition for the uniqueness of nonlinear participation factors and five other variants of participation factors, accounting for uncertain scaling factors in a mode's shape and composition. These scaling factors arise from variations in the selection of physical units or the value ranges of state variables when analyzing and controlling real-world dynamical systems. Understanding the sufficient condition of the uniqueness is therefore crucial for the correct application of participation factors in practical scenarios. Additionally, the paper explores the relationship between perturbation magnitudes in state variables and the selection of optimal scaling factors.

## 1. Introduction

In the small-signal analysis of nonlinear dynamical systems, linear participation factors (PFs) of state variables play a crucial role. These PFs are typically computed to assess the involvement of state variables in the linear modes characterized by eigenvalues of the linearized model (Garofalo et al., 2002). A linear PF is defined as the product of the corresponding elements in the right and left eigenvectors associated with an eigenvalue. This definition enables us to evaluate both the state variable's activity within the mode and its contribution to the mode itself, thus establishing a two-way connection between a state variable and a mode (Perez-arriaga et al., 1982).

In comparison, the mode shape and mode composition, two other widely used metrics defined respectively by the right and left eigenvectors of the corresponding eigenvalue, exhibit a one-way linkage and are *not* uniquely determined due to the inherent scalability of eigenvectors by any non-zero scalar (Kundur, 1993) (Sec. 12.2.2). As a common practice, the right eigenvectors (i.e. mode shapes) are often normalized, with the compositions subsequently determined based on their inverse relationship with the mode shapes (Kundur, 1993)(Eq. 12.23). Alternatively, one may normalize both mode shapes and compositions simultaneously. Importantly, even when mode shapes and compositions may not be unique due to this scaling property, linear PFs remain unique after normalization, owing to the inherent characteristics of linear systems. In the modal analysis and control of linear and nonlinear dynamical systems, the PFs of state variables with respect to a critical or selected mode serve as a pivotal tool for simplifying stability studies by focusing the system monitoring and control on a small subset of highly influential state variables (Xia, Yu & Sun et al., 2024).

Over the past two decades, researchers have introduced various types of new PFs distinct from the conventional linear PFs for stability analysis and control of dynamical systems, offering novel perspectives and applications. For instance, the concept of nonlinear PFs was introduced by leveraging the normal form theory in (Liu et al., 2006; Sanchez-Gasca et al., 2005; Shu et al., 2005), which was then applied in the design of power system controllers such as power system stabilizers to improve oscillation damping of synchronous generators under small and large disturbances. Efficient computation methods have been proposed for nonlinear PFs such as the tensor contraction-based approach in (Xia, Huang & Sun, 2024). Besides nonlinear PFs, (Abed et al., 2000) introduced the notion of probability PFs, which considers the influence of initial values and evaluates the average contribution of a mode to a state. This work explored two related variants: mode-in-state and state-in-mode probability PFs, which were subsequently examined in detail in (Hashlamoun et al., 2009). Additionally,(Hamzi & Abed, 2020) and (Iskakov, 2020) extended the concept of probability PFs to accommodate second-order nonlinearities and aspects of energy, respectively, broadening the scope of applicability. More recently, (Netto et al., 2019) adopted a formulation similar to the probability PF introduced in (Hamzi & Abed, 2020) and focused on estimating PFs using measurements within the Koopman operator-theoretic framework.

In recent years, with the increasing integration of renewable energy sources, a growing number of power electronic components have been installed in power systems. The penetration of such inverter-based resources (IBRs) significantly increases the risk of system oscillation. Therefore, participation factors and their modification indices have been applied to oscillation analysis. For example, (Yang, D., & Sun, Y., 2022) introduced a frequency-domain participation factor to identify the components with the most significant contributions and to design controllers accordingly. A similar participation factor was employed by (Yang et al., 2023) to determine the dominant device in a multi-VSC system. (Zhu et al., 2022) proposed an impedance-based participation factor to fine-tune black-box models for optimal performance, considering that many inverter-based models remain proprietary due to commercial restrictions. (Xue et al., 2023) developed a resonance participation factor using the impedance scanning method to identify the inverter-based resource with the highest contribution. Additionally, (Zhan et al., 2019) introduced loop/nodal participation factors, which consider the contribution of a loop rather than a single variable, providing a more comprehensive understanding of oscillation paths and contributing components.

The emergence of these novel types of PFs prompts a fundamental question: Do nonlinear PFs and other variants retain their uniqueness when subjected to scaling in the shape or composition of a mode? This question is crucial because, to observe and study a real-world nonlinear dynamical system, the measured or estimated values of its state variables depend on the choice of their physical units. When PFs are estimated based on a specific set of physical units, it is expected that their values may be uniquely translated to any other set of larger or smaller physical units through normalization or certain scaling factors. However, it is important to recognize that, unlike linear PFs, the PFs defined for a non-linear dynamical system, in general, cannot keep their uniqueness after the normalization of their values based on, e.g., the maximum or the sum of PFs (Dobson & Barocio, 2004; Songzhe et al., 2001). This issue can become more significant with the increase of nonlinearity of the system. For instance, in power systems, the increasing IBRs have introduced much more nonlinearities to power system dynamics. Thus, when nonlinear PFs are used to identify the highly participating devices and variables for effective control, the uniqueness of their values independent of the choice of physical units will be critical. However, the existing literature has not extensively explored this issue.

The primary objective of this paper is to identify the sufficient condition for the uniqueness of each type of PF. It is worth noting that such conditions are not straightforward and require meticulous consideration, particularly for new types of PFs, including nonlinear PFs, especially when normalization is applied in conjunction with unspecified scaling factors. Main contributions of this paper include:

1. Three types of scaling factors, namely $\xi$-factors, $\sigma$-factors, and $\theta$-factors, are introduced to represent scaling uncertainties concerning mode shapes, mode compositions, or both. It is proven that the uniqueness of most PF variants is determined by the $\theta$-factors, not individual $\xi$- or $\sigma$-factors. Specifically, the linear PF is unique if the only $\theta$-factor associated with the mode is determined (Theorem 1).

2. A sufficient condition for a nonlinear PF, as a generalization of a linear PF, to be unique to any orders of nonlinearity and combination mode is proved about all $\theta$-factors (Theorem 2), and illustrated on a toy system by two examples.

3. The relationship between the $\theta$-factors and the perturbation amplitude $\alpha$ of state variables is derived. It is shown that the perturbation amplitude $\alpha$ influences the nonlinear PF from the state variable aspect, while the scaling factor $\theta$ is viewed in terms of the mode (Remark 6).

4. It is also proved that the other five variants of PFs either share the same sufficient and necessary condition as linear PFs (with only the corresponding $\theta$-factor being determined by Corollary 1) or adhere to the same sufficient condition as nonlinear PFs (requiring the determination of all $\theta$-factors determined by Corollary 2).

The paper's primary focus lies in establishing the uniqueness condition for a nonlinear PF, as this approach simplifies the investigation of other PF variants. The paper's structure unfolds as follows: Section 2 introduces linear and nonlinear PFs; Section 3 discusses the uniqueness of the linear PF against scaling factors on eigenvectors; Section 4 presents the proof of a sufficient condition that ensures the uniqueness of nonlinear PFs of any order; Section 5 extends the proof of uniqueness conditions to encompass the remaining PF variants. Finally, Section 6 draws the conclusion.

## 2. From linear to nonlinear PFs

This section will introduce the background material, including the definitions of linear and nonlinear PFs in both non-resonant and resonant conditions.

### *2.1 Linear Participation Factor*

Consider a nonlinear dynamical system with $n$ state variables, denoted as $x_i$ ($i$ = 1, 2, ..., $n$), and a stable equilibrium located at the origin:

$$\dot{\mathbf{x}} = f(\mathbf{x}), \tag{1}$$

where state variable $\mathbf{x} \in \mathbb{R}^n$, and function $f: \mathbb{R}^n \rightarrow \mathbb{R}^n$ is assumed to be analytic. Apply the Taylor expansion at the equilibrium at the origin:

$$\dot{\mathbf{x}} = \mathbf{A}\mathbf{x} + f^{(2)}(\mathbf{x}) + f^{(3)}(\mathbf{x}) + \ldots + f^{(N)}(\mathbf{x}) + \ldots, \tag{2}$$

where $f^{(N)}(\mathbf{x})$ is the vector-valued function of all $N$-th order terms about $\mathbf{x}$ in the Taylor series (Tian et al., 2018). Assume $n$ distinct eigenvalues $\lambda_i$ with Jacobian matrix $\mathbf{A} \in \mathbb{R}^{n\times n}$, which characterize its modes. This assumption typically holds for well-designed engineering systems such as power systems operating in normal conditions (Kundur, 1993). Consider two matrices comprising the right (column) and left (row) eigenvectors of $\mathbf{A}$, respectively:

$$\mathbf{\Phi} = \begin{bmatrix} \boldsymbol{\phi}_1 & \boldsymbol{\phi}_2 & \cdots & \boldsymbol{\phi}_n \end{bmatrix}, \tag{3a}$$

$$\mathbf{\Psi} = \begin{bmatrix} \boldsymbol{\psi}_1^{\mathrm{T}} & \boldsymbol{\psi}_1^{\mathrm{T}} & \cdots & \boldsymbol{\psi}_n^{\mathrm{T}} \end{bmatrix}^{\mathrm{T}}, \tag{3b}$$

satisfying

$$\begin{cases} \mathbf{A}\boldsymbol{\phi}_i = \lambda_i \boldsymbol{\phi}_i \\ \boldsymbol{\psi}_i \mathbf{A} = \lambda_i \boldsymbol{\psi}_i \end{cases} \quad i = 1, 2, \ldots, n, \tag{3c}$$

where $\boldsymbol{\phi}_i$ and $\boldsymbol{\psi}_i$ tell the shape and composition of mode $i$ w.r.t eigenvalue $\lambda_i$, respectively (Tzounas et al., 2020).

**Definition 1:** A **linear** PF for the $k$-th state in the $i$-th mode, denoted as $p_{ki}$, is defined as the product of the $k$-th element in the $i$-th right eigenvector $\boldsymbol{\phi}_i$ and the corresponding element in the left eigenvector $\boldsymbol{\psi}_i$ of the state matrix $\mathbf{A}$ (Kundur, 1993):

$$p_{ki} \overset{\text{def}}{=} \phi_{ki} \psi_{ik}. \tag{4}$$

*Remark* 1: The linear PF, denoted as $p_{ki}$, can be interpreted as the contribution of the $i$-th mode to the $k$-th state (Kundur, 1993) or equivalently, the $k$-th state to the $i$-th mode (Hashlamoun et al., 2009) for a linear system. As demonstrated later in the paper, such interpretations are generalized and differentiated when defining various variants of PFs for a nonlinear system.

### *2.2 Nonlinear Participation Factor*

A nonlinear PF can be defined based on normal form theory (Sanchez-Gasca et al., 2005), which nonlinearly transforms the system (2) around state vector $\mathbf{x}$ into a formally linear system using a new state vector $\mathbf{z}$ by changing the coordinates in the state space (Liu et al., 2006; Shu et al., 2005). Subsequently, mode analysis can be done on this resulting linear system with the $\mathbf{z}$ state vector.

In practical applications, the normal form method is employed up to a desired order $N$ to eliminate all nonlinear terms of orders $\leq N$. Consequently, when terms of orders $> N$ are truncated, the resulting $N$-jet system becomes a linear system with respect to the new coordinates $\mathbf{z}$. While the normal form can be applied to any order, it is most commonly used in 2nd order (Sanchez-Gasca et al., 2005) or 3rd order (Amano et al., 2006; Tian et al., 2018). Below, a 2nd order nonlinear PF is introduced as an example.

First, let $\mathbf{x}$=$\mathbf{\Phi y}$ and then (2) becomes

$$\dot{y}_i = \lambda_i y_i + \sum_{p=1}^{n} \sum_{q=1}^{n} C_{pq}^{i} y_p y_q + \ldots, \tag{5}$$

where $C^i_{pq} \in \mathbb{R}^n$ denotes the coefficients of 2nd order terms after the transformation. Note that its superscript $i$ is not an exponent; rather, it represents the index of the corresponding state variable $y_i$ after the transformation (Dobson & Barocio, 2004). To eliminate 2nd order terms in (5), a nonlinear coordinate transformation **y**=**h**(**z**) is introduced (Shu et al., 2005):

$$y_i = z_i + \sum_{p=1}^{n}\sum_{q=1}^{n} h^i_{pq} z_p z_q . \tag{6}$$

Assuming there is no resonance in the system (resonance will be discussed in Section 2.3), meaning that $\lambda_p+\lambda_q-\lambda_i \neq 0$ for $\forall p$, $q$ and $i$, and if each $h$-coefficient satisfies

$$h^i_{pq} = \frac{C^i_{pq}}{\lambda_p + \lambda_q - \lambda_i}, \tag{7}$$

the resulting system in **z**-space exhibits nonlinearities of only the 3rd order or higher. A detailed proof for this transformation can be found in (Wiggins, 2003) (Chapter 19), i.e.:

$$\dot{\mathbf{z}} = \mathbf{\Lambda}\mathbf{z} + O(\|\mathbf{z}\|^3).$$

Neglecting its high-order nonlinear terms in **z**-space, the closed-form solutions in **z**, and the solutions transformed back to **y** and **x** spaces are (Liu et al., 2006)

$$z_i(t) = z_{i0} e^{\lambda_i t}, \tag{8a}$$

$$y_i(t) = z_{i0} e^{\lambda_i t} + \sum_{p=1}^{n}\sum_{q=1}^{n} h^i_{pq} z_{p0} z_{q0} e^{(\lambda_p+\lambda_q)t}, \tag{8b}$$

$$x_k(t) = \sum_{i=1}^{n} \phi_{ki} z_{i0} e^{\lambda_i t} + \sum_{i=1}^{n} \phi_{ki} [\sum_{p=1}^{n}\sum_{q=1}^{n} h^i_{pq} z_{p0} z_{q0} e^{(\lambda_p+\lambda_q)t}] \tag{8c}$$

In the case of a nonlinear system described in (2), a nonlinear PF can be defined to quantify the magnitude of mode oscillation in a state variable when only that particular state variable is perturbed. This concept is an extension of the linear PF, as discussed in Remark 1, and can be found in (Sanchez-Gasca et al., 2005) (pp. 4) and (Starrett & Fouad, 1998) (Sec. 6). An explicit expression for the 2nd order nonlinear PF is provided below.

Let initial state $\mathbf{x}_0$ have $\alpha_k$ at its $k$-th element and 0 elsewhere to represent the perturbation for the $k$-th state:

$$\mathbf{x}_0 \overset{\text{def}}{=} \begin{bmatrix} 0 & \dots & 0 & \underset{k\text{th element}}{\alpha_k} & 0 & \dots & 0 \end{bmatrix}^{\mathrm{T}}.$$

$\alpha_k$ is the perturbation amplitude for the $k$-th state variable and is commonly assumed to have a value of 1 in many papers (Sanchez-Gasca et al., 2005; Shu et al., 2005). When substituting it into (6), the initial state $z_{i0}$ is typically approximated by (Shu et al., 2005) (pp. 4).

$$z_{i0} = \alpha_k \psi_{ik} - \alpha_k^2 \sum_{p=1}^{n} \sum_{q=p}^{n} h_{pq}^{i} \psi_{pk} \psi_{qk}. \tag{9}$$

In (9), the index $q$ starts from $p$, which is a common practice in the calculation of nonlinear PFs. A detailed discussion concerning this index can be found in (Sanchez-Gasca et al., 2005) (Sec. II-A). Plugging (9) into (8c), the closed-form solution is obtained:

$$x_k(t) = \sum_{i=1}^{n} p_{2ki} e^{\lambda_i t} + \sum_{p=1}^{n} \sum_{q=p}^{n} p_{2kpq} e^{(\lambda_p + \lambda_q)t}, \tag{10}$$

$$p_{2ki} = \phi_{ki}(\alpha_k \psi_{ik} + \psi_{2ikk}) = \alpha_k p_{ki} + \alpha_k^2 p_{2kiNL}, \tag{11a}$$

$$p_{2kpq} = \phi_{2kpq}(\psi_{pk} + \psi_{2pkk})(\psi_{qk} + \psi_{2qkk}), \tag{11b}$$

where

$$\psi_{2mkk} = -\alpha_k^2 \sum_{p=1}^{n} \sum_{q=p}^{n} h_{pq}^{m} \psi_{pk} \psi_{qk} \quad , \quad \phi_{2kpq} = \sum_{i=1}^{n} h_{pq}^{i} \phi_{ki}.$$

*Remark* 2: The two equations in (11) provide formulas for two variants of PFs that account for 2nd order nonlinearities. In (11a), $p_{2ki}$ is defined as the 2nd order nonlinear PF of the $k$-th state variable in linear mode $i$, which equals the linear PF $p_{ki}$ multiplied by the perturbation amplitude $\alpha_k$, along with an additional correction term $\alpha_k^2 p_{2kiNL}$. Regarding $p_{2kpq}$ in (11b), it represents the nonlinear PF of the $k$-th state variable in a combination mode characterized by two linear modes $\lambda_p + \lambda_q$ (Amano et al., 2006). It's worth noting that although such a mode is named as a 2nd order nonlinear mode in some literature, such as (Liu et al., 2006), this paper follows the task force report (Sanchez-Gasca et al., 2005) and terms it as the combination mode.

When $\alpha_k$ =1, or equivalently, $\mathbf{x}_0 = \mathbf{e}_k$, the first term in Equation (11a) becomes identical to the linear PF $p_k$. Some researchers (Liu et al., 2006) prefer to retain this unit perturbation to preserve this consistency property. This property is also maintained in report (Sanchez-Gasca et al., 2005) and is widely adopted in the literature. Although (Shu et al., 2005) identified this issue and introduced $\alpha$ as a relaxation parameter, they merely suggested selecting a suitable value without providing a detailed discussion. In (Xia, T., & Sun, K., 2022), the primary focus was on establishing a connection between linear and nonlinear participation factors rather than on their uniqueness. In practical systems, considering unit and base values in a per-unit system, it is often more prudent to keep $\alpha_k$ as a variable rather than fixing its value 1 during formula derivation. This approach facilitates a better understanding of the scaling factor's impact, as demonstrated by Example 2 in Section 4.

### *2.3 On Resonance*

A first-order resonance, often called a strong resonance, occurs when the state matrix A has two identical eigenvalues (Dobson & Barocio, 2005). (11) remains valid even if the Jordan canonical form is employed for non-diagonalizable A, as described in (Sanchez-Gasca et al., 2005) (Eq. 4), based on a generalization of Poincare's theorem (Arnold, 1988) (Sec. 23C).

In well-designed real-life systems, it's not common for the eigenvalues to be exactly equal, and therefore, strong resonance is not a common occurrence. However, near resonance can arise when two eigenvalues are very close to each other, and Detailed studies can be found in (Dobson et al., 2001).

A 2$^{nd}$ order resonance occurs when $\lambda_p+\lambda_q-\lambda_i = 0$, $\exists$ $p$, $q$ and $i$. Additionally, real-life engineering systems, such as power systems, can have zero eigenvalues, which constitute a special type of 2$^{nd}$ order resonance (Samovol, 2004) (Theorem 3). Unfortunately, the definition of the nonlinear PF under resonant conditions is not found in existing literature. Nevertheless, the response of a system with resonance can still be approximated using (Wang & Huang, 2017) (Eq. 19):

$$x_k(t) = \sum_{i=1}^{n} \phi_{ki} z_{i0} e^{\lambda_i t} + \underbrace{\sum_{i=1}^{n} \phi_{ki} \sum_{p=1}^{n} \sum_{q=1}^{n} h_{pq}^{i} z_{p0} z_{q0} e^{(\lambda_p+\lambda_q)t}}_{\lambda_p+\lambda_q \neq \lambda_i} + \underbrace{\sum_{i=1}^{n} \phi_{ki} \sum_{p=1}^{n} \sum_{q=1}^{n} C_{pq}^{i} z_{p0} z_{q0} (1+t) e^{\lambda_i t}}_{\lambda_p+\lambda_q = \lambda_i} .$$

which introduces a third term that grows with time compared to (8c). It will become evident later that even when considering resonance or near resonance, the conclusions regarding nonlinear PFs in this paper remain valid based on (20). This is because the factor $\lambda_p+\lambda_q-\lambda_i$ or $1+t$ does not affect the scaling of eigenvectors. Although we only demonstrate the case of 2$^{nd}$ order resonance here, scenarios with higher-order resonance lead to similar conclusions.

## 3. From Linear Systems to Nonlinear Systems

This section establishes the uniqueness of a linear PF against scaling uncertainties in mode shape and mode composition by introducing three *scaling factors*: $\xi$-factors, $\sigma$-factors and $\theta$-factors, which respectively scale mode shapes, mode compositions, and both.

### *3.1Scaling Factors*

If $\boldsymbol{\phi}_i$ is a right eigenvector (mode shape) of $\lambda_i$, it remains so after being scaled by any non-zero scalar (Kundur, 1993) (Sec. 12.2.2). Without loss of generality, we define unique mode shapes and mode compositions, each with a *unit* norm:

$$\hat{\boldsymbol{\phi}}_i = \frac{\boldsymbol{\phi}_i}{\|\boldsymbol{\phi}_i\|} \quad , \quad \hat{\boldsymbol{\psi}}_i = \frac{\boldsymbol{\psi}_i}{\|\boldsymbol{\psi}_i\|}. \tag{12}$$

Let $\hat{\boldsymbol{\phi}}_i$ and $\hat{\boldsymbol{\psi}}_i$ be the $i$-th right (column) and left (row) eigenvectors, each with a unit norm (e.g., a unity 2-norm). There exist unique **scaling factors** $\sigma_i$ and $\xi_i \in \mathbb{C}$ such that, for any left and right eigenvectors $\boldsymbol{\Phi}$ and $\boldsymbol{\Psi}$ in (3), the following holds:

$$\boldsymbol{\Phi} = \left[ \sigma_1 \hat{\boldsymbol{\phi}}_1 \quad \sigma_2 \hat{\boldsymbol{\phi}}_2 \quad \dots \quad \sigma_n \hat{\boldsymbol{\phi}}_n \right], \tag{13a}$$

$$\boldsymbol{\Psi} = \left[ \xi_1 \hat{\boldsymbol{\psi}}_1^{\mathrm{T}} \quad \xi_2 \hat{\boldsymbol{\psi}}_2^{\mathrm{T}} \quad \dots \quad \xi_n \hat{\boldsymbol{\psi}}_n^{\mathrm{T}} \right]^{\mathrm{T}}, \tag{13b}$$

Thus, the products based on those $\boldsymbol{\Phi}$ and $\boldsymbol{\Psi}$ given by

$$\theta_i \overset{\text{def}}{=} \boldsymbol{\psi}_i \boldsymbol{\phi}_i = \xi_i \hat{\boldsymbol{\psi}}_i \sigma_i \hat{\boldsymbol{\phi}}_i = \xi_i \sigma_i \left( \cos \delta_i \right), \tag{13c}$$

where $\delta_i$ represents the angle between the mode shape $\boldsymbol{\psi}_i$ and mode composition $\boldsymbol{\phi}_i$. Throughout the rest of this paper, the sets of $\sigma_i$, $\xi_i$ and $\theta_i$ ($i$=1, …, $n$) are referred to as $\xi$-factors, $\sigma$-factors and $\theta$-factors, respectively. If $\theta_i = 1$ for any $i$, it implies that $\boldsymbol{\Psi} = \boldsymbol{\Phi}^{-1}$. From (13c), as $\delta_i$ is a constant for a particular system, the scaling factors $\sigma_i$ and $\xi_i$ uniquely determine the value of $\theta_i$.

*Remark* 3: By introducing the scaling factors, any other mode shape and mode composition matrices can be expressed using $\hat{\boldsymbol{\phi}}_i$ and $\hat{\boldsymbol{\psi}}_i$ with scaling factors $\sigma_i$ and $\xi_i$. Without specified notation, the norm in the following discussion refers to the 2-norm, as paper (Kundur, 1993; Liu et al., 2006; Sanchez-Gasca et al., 2005). In fact, extending it to the $p$-norm does not affect the conclusions in this paper. Additionally, the mode shapes and mode compositions for different modes are orthogonal, resulting in their inner product being equal to zero (Kundur, 1993) (Eq. 12.21).

### *3.2 Uniqueness of the Linear PF*

Based on the definitions of linear PF in (4) and scaling factor in (13), we have the following expression:

$$\mathbf{P} = \boldsymbol{\Phi} \circ \boldsymbol{\Psi}^{\mathrm{T}} = \left[ \theta_1 \left( \frac{\hat{\boldsymbol{\phi}}_1 \circ \hat{\boldsymbol{\psi}}_1^{\mathrm{T}}}{\cos \delta_1} \right) \quad \theta_2 \left( \frac{\hat{\boldsymbol{\phi}}_2 \circ \hat{\boldsymbol{\psi}}_2^{\mathrm{T}}}{\cos \delta_2} \right) \quad \dots \quad \theta_n \left( \frac{\hat{\boldsymbol{\phi}}_n \circ \hat{\boldsymbol{\psi}}_n^{\mathrm{T}}}{\cos \delta_n} \right) \right], \tag{14}$$

where "∘" denotes the Hadamard product, which represents element-wise multiplication. The $i$-th column of matrix $\mathbf{P}$ contains the linear PFs of all state variables associated with mode $i$ for a given $\theta_i$. Consequently, the following sufficient and necessary condition for unique linear PFs with each mode can be derived:

**Theorem 1**: Providing a scaling factor $\theta_i \in \mathbb{C}$ defined in (13), with $i \in \{1, 2, ..., n\}$, the linear PFs of all state variables associated with mode $i$, denoted as $p_{ki}$ for $k \in \{1, 2, ..., n\}$, are unique *if and only if the corresponding* $\theta_i$ *is* unique.

Theorem 1 indicates that the linear PFs associated with mode *i* are unique if and only if $\theta_i$ is determined. This theorem highlights a crucial property of linear PFs within each mode *i*: the linear PF for each state variable remains constant regardless of changes in $\sigma_i$ or $\xi_i$ once their product $\theta_i$ is fixed. Consequently, the vector of linear PFs for all state variables within mode *i* becomes unique after normalization.

**Example 1:** Consider a linear system shown by

$$\begin{bmatrix} \dot{x}_1 \\ \dot{x}_2 \\ \dot{x}_3 \\ \dot{x}_4 \end{bmatrix} = \mathbf{A}\mathbf{x} = \begin{bmatrix} 0 & 0 & 1 & 0 \\ 0 & 0 & 0 & 1 \\ -20 & 20 & -1 & 0 \\ 5 & -5 & 0 & -1 \end{bmatrix} \begin{bmatrix} x_1 \\ x_2 \\ x_3 \\ x_4 \end{bmatrix},$$

where **A** is the state matrix with four eigenvalues: $\lambda_1 = -0.50+4.97j$, $\lambda_2 = -0.50-4.97j$, $\lambda_3 = 0$ and $\lambda_4 = -1.00$. Three cases are:

I) $\sigma_i = 1$ and $\xi_i$=1 for any *i*, i.e., normalizing mode shapes and compositions, respectively, to have a unity norm.

II) $\sigma_i = 1$ and $\boldsymbol{\psi}_i\boldsymbol{\phi}_i = 1$.

III) $\xi_i = 1$ and $\boldsymbol{\psi}_i\boldsymbol{\phi}_i = 1$.

The scaling factors for the four eigenvalues are displayed in Table 1. This study primarily focuses on the oscillation mode with $\lambda_1 = -0.50+4.97j$ to conserve space. With the corresponding scaling factors $\sigma_1$, $\xi_1$, and $\theta_1$, the mode shapes and mode compositions for Case I, II, and III are provided in Table 2. It's worth noting that the directions of mode shapes or mode compositions remain the same across all three cases, while the amplitudes differ due to the scaling factors.

**Table 1.** The scaling factors for four eigenvalues in Case I, II and III

| # | $\sigma_1$ | $\sigma_2$ | $\sigma_3$ | $\sigma_4$ | $\xi_1$ | $\xi_2$ | $\xi_3$ | $\xi_4$ | $\theta_1$ | $\theta_2$ | $\theta_3$ | $\theta_4$ |
|---|---|---|---|---|---|---|---|---|---|---|---|---|
| Case I | 1.000 | 1.000 | 1.000 | 1.000 | 1.000 | 1.000 | 1.000 | 1.000 | 0.138 | 0.138 | 0.250 | 0.250 |
| Case II | 1.000 | 1.000 | 1.000 | 1.000 | 7.236 | 7.236 | 4.000 | 4.000 | 1.000 | 1.000 | 1.000 | 1.000 |
| Case III | 7.236 | 7.236 | 4.000 | 4.000 | 1.000 | 1.000 | 1.000 | 1.000 | 1.000 | 1.000 | 1.000 | 1.000 |

**Table 2.** The mode shapes and compositions of $\lambda_1$ in Case I, II and III

| | Case I | Case II | Case III | Legends |
|---|---|---|---|---|
| Mode Shapes (Right Eigenvectors) | $\boldsymbol{\phi}_1^{\mathrm{I}} = \sigma_1^{\mathrm{I}} \hat{\boldsymbol{\phi}}$ | $\boldsymbol{\phi}_1^{\mathrm{II}} = \sigma_1^{\mathrm{II}} \hat{\boldsymbol{\phi}}$ | $\boldsymbol{\phi}_1^{\mathrm{III}} = \sigma_1^{\mathrm{III}} \hat{\boldsymbol{\phi}}_1$ | State index: $x_1$, $x_2$, $x_3$, $x_4$ |
| Mode Compositions (Left Eigenvectors) | $\boldsymbol{\psi}_1^{\mathrm{I}} = \xi_1^{\mathrm{I}} \hat{\boldsymbol{\psi}}_1$ | $\boldsymbol{\psi}_1^{\mathrm{II}} = \xi_1^{\mathrm{II}} \hat{\boldsymbol{\psi}}_1$ | $\boldsymbol{\psi}_1^{\mathrm{III}} = \xi_1^{\mathrm{III}} \hat{\boldsymbol{\psi}}_1$ | Amplitude: 0.6, 0.4, 0.2 |

Applying (14), the linear PFs for the mode are

$$p_{k1}^{\mathrm{I}} = \begin{bmatrix} 0.056 \\ 0.014 \\ 0.056 \\ 0.014 \end{bmatrix}, \quad p_{k1}^{\mathrm{II}} = p_{k1}^{\mathrm{III}} = \begin{bmatrix} 0.402 \\ 0.101 \\ 0.402 \\ 0.101 \end{bmatrix} \xrightarrow{\text{Normalization}} \begin{bmatrix} 1.000 \\ 0.250 \\ 1.000 \\ 0.250 \end{bmatrix}.$$

Notice that $p_{k1}^{\mathrm{I}}/p_{k1}^{\mathrm{II}} = \theta_1^{\mathrm{I}}/\theta_1^{\mathrm{II}}$ for any $k$ before the normalization. In the case of a linear system, the PFs for mode $i$ are directly proportional to $\theta_i$. Consequently, after normalization, all PFs are equal to the same vector. Figure 1 illustrates the responses of the linear system under a specific perturbation. The left figure represents the response of the $k$-th state variable when only the $k$-th state variable is perturbed with $\alpha_k$, where $\alpha_k = \theta_1^{\mathrm{II}} = 1$. Each response consists of four eigenvalue components ($\lambda_1$ to $\lambda_4$) based on linear system theory; the components for $\lambda_1$ are shown in the right figure, where the amplitude of each oscillation (envelope) is just the PF for each state in Case II based on the physical meaning of PF; Case III yields the same results due to the identical $\theta$-factor. For this linear system, the response with $\alpha_k = \theta_1^{\mathrm{I}} = 0.138$ closely resembles the left figure and the component for $\lambda_1$ in Case I is shown in the middle figure.

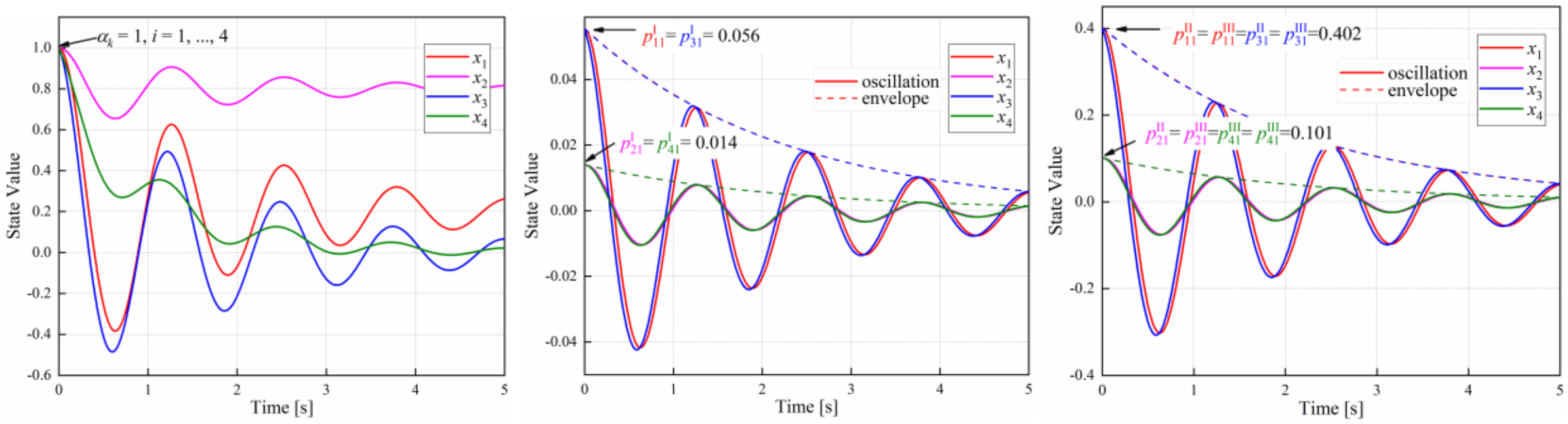

**Figure 1.** The responses of the linear system

(Left: the responses for each state variable when perturbation amplitude $\alpha_k$ = $\theta_1^{\mathrm{II}}$=1; Middle: the component for oscillation mode $\lambda_1$ = −0.50+4.97$j$ in Case I; Right: the component for oscillation mode $\lambda_1$ = −0.50+4.97$j$ in Case

## 4. Uniqueness of The Nonlinear PFs

In this section, we establish a sufficient condition for the uniqueness of a nonlinear PF of any nonlinearity order with a linear or combination mode (as discussed in Remark 2) in the presence of scaling uncertainties in eigenvectors. We begin by deriving a general expression for a nonlinear PF and subsequently provide a detailed proof of its uniqueness. To illustrate these concepts, an example is presented, and further, the relationships between scaling factors and perturbation amplitudes are explored.

### *4.1 Normal Form Transformation*

Notation: To obtain a theorem covering any order of nonlinearities for both linear and combination modes, it is essential to clarify the *orders* of a nonlinear PF and a mode. In the following content, $N \in \mathbb{Z}^+$ is used to represent the nonlinearity order, corresponding to the order of the highest nonlinearity considered in the Taylor series. Additionally, $M \in \mathbb{Z}^+$ is employed to denote the combination order of the combination mode, where $M = 1$ signifies a linear mode. For instance, in a 2nd order nonlinear PF, as depicted in (11), $N$ is fixed at 2 to truncate terms with nonlinearities of orders greater than 2, while $M = 1$ for (11a) and $M = 2$ for (11b). It is worth noting that, due to the utilization of the normal form method, $M$ is constrained by the order of the Taylor series, resulting in $M \leq N$.

The Taylor expansion of (1) up to an infinite order is represented as follows:

$$\dot{x}_k = \sum_{i=1}^{n} a_{ki} x_i + \sum_{p=1}^{n} \sum_{q=1}^{n} a_{k,pq} x_p x_q + \ldots + \sum_{r=1}^{n} \ldots \sum_{v=1}^{n} a_{k,\underbrace{r \ldots v}_{N}} x_r \ldots x_v + \cdots, \tag{15}$$

where $x_k$ denotes the $k$-th state variable, $a_{ki}$ represents the element in the $k$-th row and $i$-th column of state matrix $\mathbf{A}$, $a_{k,pq}$ is the $p$-th row and $q$-th column element in $k$-th Hessian matrix, $a_{k,r \ldots v}$ is the coefficient of $N$-th order Taylor series term.

Similar to (8c), a closed-form expression in **x**-space up to an infinity order is given by (Huang et al., 2009)

$$x_k = \sum_{i=1}^{n} \phi_{ki} z_i + \sum_{i=1}^{n} \phi_{ki} \sum_{p=1}^{n} \sum_{q=1}^{n} h_{pq}^{i} z_p z_q + \ldots + \sum_{i=1}^{n} \phi_{ki} \sum_{r=1}^{n} \ldots \sum_{v=1}^{n} h_{\underbrace{rs\ldots v}_{N}}^{i} z_r z_s \ldots z_v + \ldots, \tag{16a}$$

where

$$z_i = z_{i0} e^{\lambda_i t}, \tag{16b}$$

$$h_{\underbrace{rs\ldots v}_{N}}^{i} \approx \frac{\sum_{j=1}^{n} \sum_{\alpha=1}^{n} \ldots \sum_{\gamma=1}^{n} \psi_{ij} a_{j,\alpha\beta\ldots\gamma} \phi_{\alpha r} \phi_{\beta s} \ldots \phi_{\gamma v}}{\lambda_r + \lambda_s + \ldots + \lambda_v - \lambda_i}. \tag{16c}$$

Each $h$-coefficient has an approximate formula in (16c) if the $N$-th order nonlinearity is significantly more dominant than any lower-order nonlinearities. Let $N = 2$, the (16) will be downgraded to 2nd order normal form where (16a) corresponds to (8c), (16b) is identical to (8a) and (16c) becomes (7).

### ***4.2Nonlinear PF of any Nonlinearity and Combination Orders***

Since the normal form expression for any order has been derived in (16), the corresponding nonlinear PF will be derived in this part. We continue following the definition of nonlinear PF in (11) and retain the perturbation amplitude $\alpha_k$. Similar to (9), the initial state in **z**-space is given by:

$$z_{i0} = \alpha_k \psi_{ik} - \ldots - \alpha_k^{N} \sum_{r=1}^{n} \ldots \sum_{w}^{n} \sum_{v=w}^{n} h_{\underbrace{rs\ldots v}_{N}}^{i} \psi_{rk} \psi_{sk} \ldots \psi_{vk} - \ldots,$$

Notice that the index $r$ starts from 1 while the index $v$ starts from the index $w$, just as $q$ starts from $p$ in (9). For simplicity, this initial value expression can be rewritten as

$$\mu_{ik} = z_{i0} \big|_{\mathbf{x}_0 = \mathbf{e}_k}. \tag{17}$$

Similar to (11), the nonlinear PF with a linear mode is

$$p_{ki} = \phi_{ki} \mu_{ik}, \tag{18a}$$

and the nonlinear PF for an $M$-th ($M \leq N$) order combination mode involving $M$ modes with indices $r$, $s$, …, $u$ is

$$p_{k,\underbrace{rs\ldots u}_{M}} = \sum_{i=1}^{n} \phi_{ki} h_{rs\ldots u}^{i} \underbrace{\mu_{rk} \mu_{sk} \ldots \mu_{uk}}_{M}. \tag{18b}$$

If $h_i^i$ is set to 1 for $i = 1$, the nonlinear PF in (18a) can be regarded as a particular case of (18b) with a combination mode order when M = 1.

Note that (16a) contains an infinite number of terms, allowing (18) to define a nonlinear PF considering nonlinearities of any order. In practice, calculating a nonlinear PF is typically performed up to a desired order *N*, with all terms of orders greater than *N* truncated. For instance, when $N = 2$, (18) yields the same 2nd order nonlinear PF as defined in (11), while a 3rd order example is provided in (Tian et al., 2018) and (Huang et al., 2009).

### *4.3 Uniqueness of a Nonlinear PF*

There is the following theorem on the uniqueness of nonlinear PFs up to *N*-th order.

**Theorem 2:** For a scaling factor $\theta_i \in \mathbb{C}$ defined in (13), $i \in \{1, 2, ..., n\}$, the nonlinear PFs represented by $p_{k,rs\ldots u}$ with $k \in \{1, 2, ..., n\}$, for all state variables associated with a single linear or combination mode constructed by *M* modes *r, s…,u*, are unique if *all θ-factors are unique*.

*Proof:*

In the following proof, any variable with a hat (^) signifies its *irrelevance* from the scaling factors $\xi_i$, $\sigma_i$ or $\theta_i$. According to (13), the mode shape and composition with respect to mode *i* are as follows:

$$\phi_{ki} = \sigma_i \hat{\phi}_{ki} \text{ and } \psi_{ik} = \xi_i \hat{\psi}_{ik}. \tag{19}$$

For an *M*-th order combination mode (or a linear mode for $M = 1$), substitute (19) into (18b), and the nonlinear PF becomes

$$\underset{M}{p_{k,rs\ldots u}} = \sum_{i=1}^{n} \sigma_i \hat{\phi}_{ki} h^i_{rs\ldots u} \underbrace{\mu_{rk}\mu_{sk}\cdots\mu_{uk}}_{M}, \tag{20a}$$

$$\mu_{lk} = \alpha_k \xi_l \hat{\psi}_{lk} - \ldots - \alpha_k^N \sum_{r=1}^{n} \cdots \sum_{v=u}^{n} \left(\xi_r \xi_s \ldots \xi_v\right) \underset{N}{h^l_{rs\ldots v}} \hat{\psi}_{rk} \hat{\psi}_{sk} \ldots \hat{\psi}_{vk} - \ldots , \tag{20b}$$

$\forall\, l \in \{r, s, \ldots, u\}$. Note that, in (20b)

$$\underset{N}{h^l_{rs\ldots v}} = \left(\xi_l \sigma_r \sigma_s \ldots \sigma_v\right) \underbrace{\frac{\sum_{j=1}^{n} \sum_{\alpha=1}^{n} \cdots \sum_{\gamma=1}^{n} \hat{\psi}_{lj} a_{j,\alpha\beta\ldots\gamma} \hat{\phi}_{\alpha r} \hat{\phi}_{\beta s} \ldots \hat{\phi}_{\gamma v}}{\lambda_r + \lambda_s + \ldots + \lambda_v - \lambda_l}}_{\hat{h}^l_{rs\ldots v}}. \tag{20c}$$

By substituting (20c) into (20b), all *σ*-factors and *ξ*-factors can be replaced by their products *θ*-factors except for $\xi_l$, as demonstrated below

$$\mu_{lk} = \xi_l \alpha_k \hat{\psi}_{lk} - \ldots - \alpha_k^N \sum_{r=1}^{n} \ldots \sum_{v=w}^{n} \left(\xi_r \xi_s \ldots \xi_v\right)\left(\xi_l \sigma_r \sigma_s \ldots \sigma_v\right) \underset{N}{\hat{h}^l_{rs\ldots v}} \hat{\psi}_{rk} \hat{\psi}_{sk} \ldots \hat{\psi}_{vk} - \ldots$$

$$= \xi_l \left( \alpha_k \hat{\psi}_{lk} - \ldots - \alpha_k^N \sum_{r=1}^{n} \ldots \sum_{v=w}^{n} \underbrace{\frac{\theta_r \theta_s \ldots \theta_v}{\cos\delta_r \cos\delta_s \ldots \cos\delta_v}}_{N} \hat{h}^l_{rs\ldots v} \underbrace{\hat{\psi}_{rk} \hat{\psi}_{sk} \ldots \hat{\psi}_{vk}}_{N} - \ldots \right).$$

Similarly, as in the case of (20c), the coefficient $h^i_{rs\ldots u}$ in (20a) can be written as

$$\underset{M}{h^i_{rs\ldots u}} = \left(\xi_i \sigma_r \sigma_s \ldots \sigma_u\right) \hat{h}^i_{rs\ldots u}.$$

Therefore, substituting it into (20a) allows us to eliminate all scaling factors $\sigma_i$ and $\xi_i$. Consequently, (20) simplifies to

$$\underset{M}{p_{k,rs\ldots u}} = \sum_{i=1}^{n} \frac{\theta_i}{\cos\delta_i} \hat{\phi}_{ki} \hat{h}^i_{rs\ldots u} \underbrace{\mu_{rk} \mu_{sk} \ldots \mu_{uk}}_{M}, \tag{21a}$$

$$\mu_{lk} = \frac{\theta_l}{\cos\delta_l} \left( \alpha_k \hat{\psi}_{lk} - \ldots - \alpha_k^N \sum_{r=1}^{n} \ldots \sum_{v=w}^{n} \underbrace{\frac{\theta_r \theta_s \ldots \theta_v}{\cos\delta_r \cos\delta_s \ldots \cos\delta_v}}_{N} \hat{h}^l_{rs\ldots v} \hat{\psi}_{rk} \hat{\psi}_{sk} \ldots \hat{\psi}_{vk} - \ldots \right). \tag{21b}$$

Notably, the nonlinear PF is independent of both $\sigma$-factors or $\xi$-factors; rather, it relies on the determination of every $\theta_i$, $i \in \{1, 2, ..., n\}$, or in other words, the values of all $\theta$-factors. Consequently, the uniqueness of the nonlinear PF is established if and only if all $\theta$-factors are determined.

□

*Remark 4*: Theorem 2 establishes that the uniqueness of nonlinear PFs depends on the determination of all $\theta$-factors. This differs from the case of linear PFs in Theorem 1, where only the corresponding $\theta_i$ is necessary for uniqueness.

*Remark* 5: Unlike Theorem 1 for linear PFs, the condition in Theorem 2 is sufficient but not necessary. This implies that a unique set of $p_{k,rs\ldots u}$ ($k$ = 1, 2, ..., $n$) might correspond to multiple sets of $\theta$-factors.

In fact, all possible $\theta$-factors that yield a unique set of nonlinear PFs can be determined by solving the following equations:

$$\begin{cases} g_1(\theta_1, \theta_2, \ldots, \theta_n) = 0 \\ \quad\vdots \\ g_n(\theta_1, \theta_2, \ldots, \theta_n) = 0, \end{cases} \tag{22a}$$

$$g_i(\theta_1, \theta_2, \ldots, \theta_n) = \underset{M}{p_{k,rs\ldots u}} - \sum_{i=1}^{n} \frac{\theta_i}{\cos\delta_i} \hat{\phi}_{ki} \hat{h}^i_{rs\ldots u} \mu_{rk} \mu_{uk} \ldots \mu_{uk}. \tag{22b}$$

The number of roots of these equations corresponds to their Bézout number, denoted as $\prod d_i$ for all $i$ (Dreesen et al., 2012) (Th. 1), where $d_i$ represents the degree of $g_i$, which is also the order of the nonlinear PF. It is evident that as the order of nonlinear PFs increases, the solutions for $\theta$-factors will be non-unique.

It is important to note that the Bezout number accounts for complex roots. The precise count of all real roots can be determined by solving (22) (Benallou et al., 1983), which is known to be an NP-complete problem (Courtois et al., 2000). In fact, selecting one set of $\theta$-factors by Theorem 2 is sufficient to ensure the uniqueness of the calculated PFs while finding all possible sets of $\theta$-factors is rarely needed in practical applications and is not the focus of this paper.

*Remark* 6: The influence of the scaling factors $\theta_i$ is similar (though not equivalent) to adjusting the perturbation amplitude. (21b) can be expressed as:

$$\mu_{lk} = \left(\frac{\alpha_k \theta_l}{\cos \delta_l}\right)\hat{\psi}_{lk} - \ldots - \sum_{r=1}^{n}\cdots\sum_{v=w}^{n}\underbrace{\prod_{\gamma=r}^{v}\left(\frac{\alpha_k \theta_\gamma}{\cos \delta_\gamma}\right)}_{N}\hat{h}^{l}_{rs\ldots v}\hat{\psi}_{rk}\hat{\psi}_{sk}\ldots\hat{\psi}_{vk} - \ldots .$$

Therefore, modifying the perturbation amplitude $\alpha_k$ and carefully designing $\theta_i$ may either cancel each other out or produce the same effect on nonlinear PFs. Moreover, selecting different scaling factors $\theta_i$ allows for the amplification or reduction of the contribution of a specific mode $i$. Unfortunately, the issue of determining a reasonable perturbation amplitude $\alpha_k$ remains unsolved, and it is typically based on empirical knowledge (Liu et al., 2006) or set to a unit value for simplification (Perez-arriaga et al., 1982). In contrast to the perturbation amplitude $\alpha_k$, the scaling factor $\theta_i$ offers an additional dimension for adjusting nonlinear PFs. It's worth noting that the perturbation amplitude remains the same for the $k$-th state variable (i.e., $\alpha_k$), while the scaling factor remains consistent for the i-th mode (i.e., $\theta_i$).

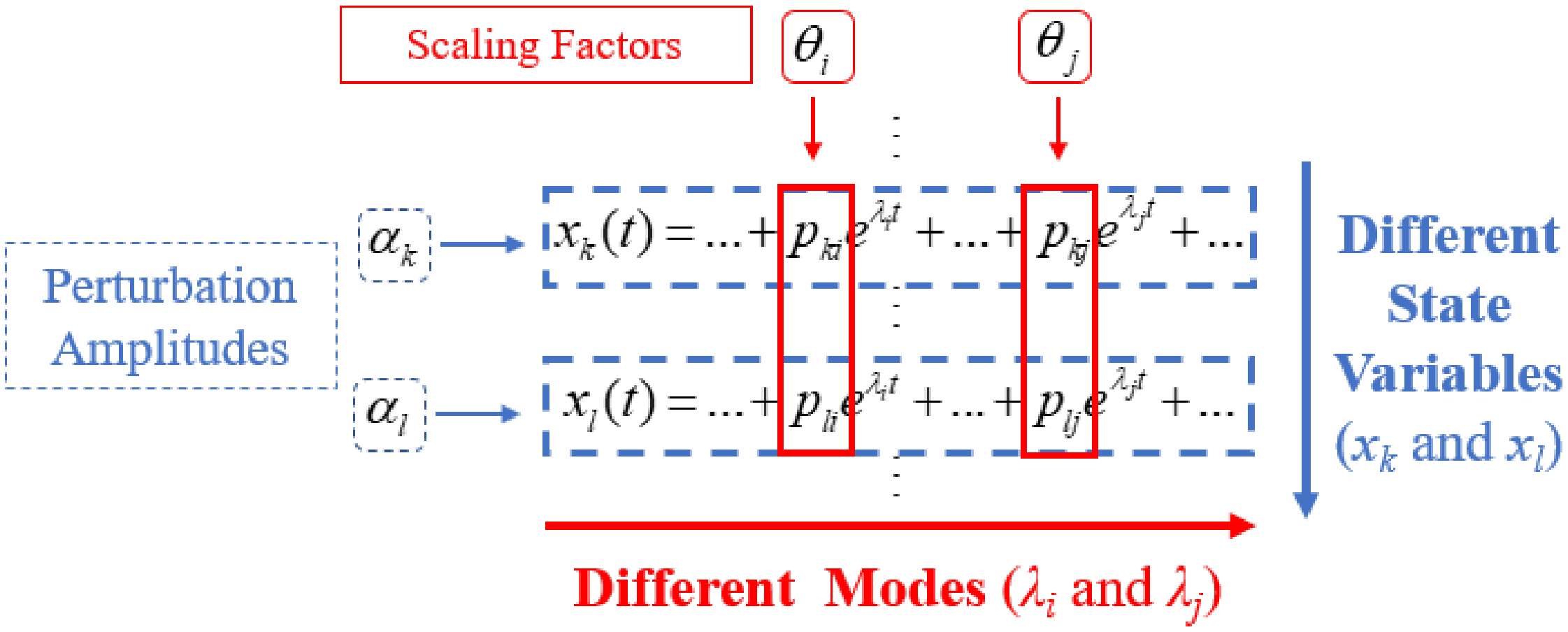

**Figure 2.** The relationship between perturbation amplitudes and scaling factors

The study of participation factors examines the relationship between state variables and modes. In (Hashlamoun et al., 2009) (Sec. 2), the authors delve into the details of state-in-mode and mode-in-state PFs, highlighting a degree of symmetry between state variables and modes. Similarly, considering that the perturbation amplitude is chosen with respect to state variables, it is logical to introduce another factor that accounts for modes, which is the scaling factor in our analysis. Figure 2 provides an illustration of this concept. Mathematically, there is no inherent reason to believe that some variables are more important than others, hence the common practice of setting $\alpha_k = \alpha$ for all $k$. Therefore, in this paper, we propose $\theta_i = \theta$. In Case II and III in Example 1, $\theta_i = 1$ for any $i$.

**Example 2:** Consider a nonlinear system with the same linear part as Example 1:

$$\dot{\mathbf{x}} = \mathbf{A}\mathbf{x} + \mathbf{g}(\mathbf{x}), \quad \mathbf{g}(\mathbf{x}) = \begin{bmatrix} 0 & 0 & -2x_1x_3 & 0 \end{bmatrix}^{\mathrm{T}}.$$

Still test three cases in Table 1. Following the normalization process (the initial values are depicted in Figure 4), we focus on the same linear mode $\lambda_1 = -0.50+4.97j$. Notably, in Case I, the PFs differ from those in Cases II and III due to variations in the $\theta$-factors:

$$p_{2,k1}^{\mathrm{I}} = \begin{bmatrix} 0.994 \\ 0.249 \\ 1.000 \\ 0.253 \end{bmatrix} \neq p_{2,k1}^{\mathrm{II}} = p_{2,k3}^{\mathrm{III}} = \begin{bmatrix} 0.865 \\ 0.222 \\ 1.000 \\ 0.353 \end{bmatrix}.$$

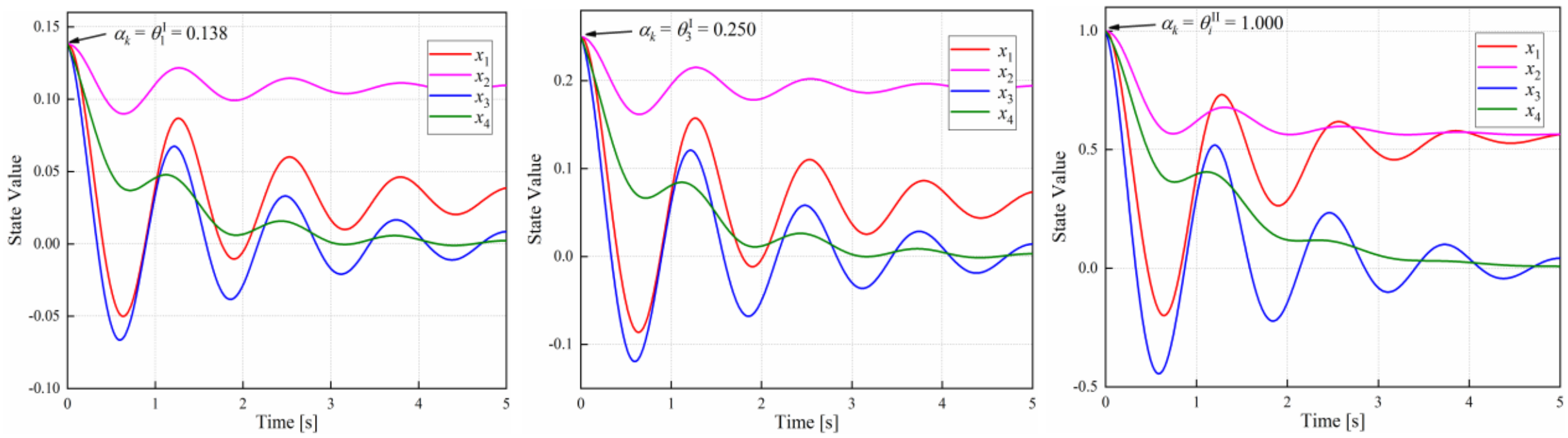

**Figure 3.** The responses of the nonlinear system with different perturbation amplitudes (Left: $\alpha_k = \theta_1^{\mathrm{I}} = 0.138$; Middle: $\alpha_k = \theta_3^{\mathrm{I}} = 0.250$; Right: $\alpha_k = \theta_i^{\mathrm{II}} = 1.000$)

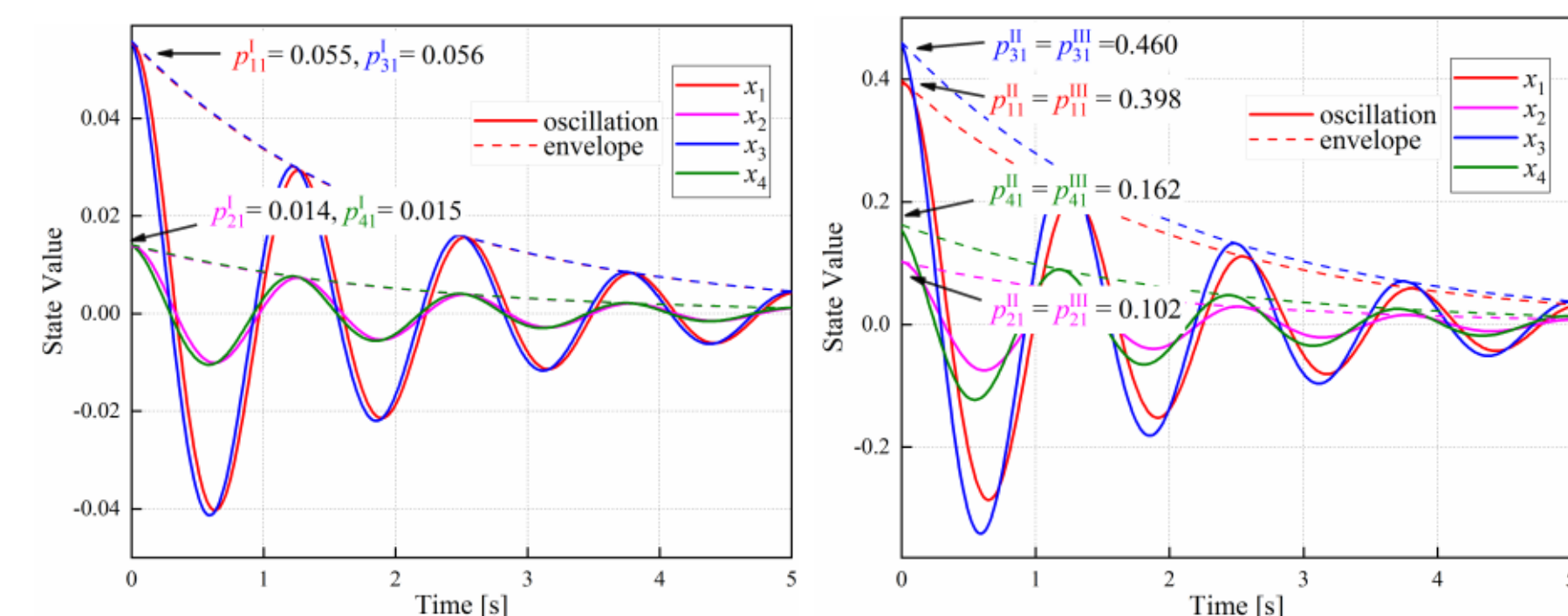

**Figure 4.** The reconstructed responses for mode $\lambda_1$ in Case I (left) and Case II or III (right)

The examination focuses on the system's responses under varying perturbation amplitudes to clarify the observed differences, as depicted in Figure 3. In Cases II or III (where the results are the same), the scaling factor $\theta_i$ remains consistent across all four modes ($\lambda_1$ to $\lambda_4$). Consequently, following the definition, the perturbation amplitude is set as $\theta_i^{\mathrm{II}}=1$ to illustrate the nonlinear PF. However, in Case I, the scaling factor $\theta_i$ varies among the four modes, rendering it impractical to represent them under a single type of disturbance. In this scenario, two distinct perturbation amplitudes are employed, and their responses resemble those of the linear system shown in Figure 1.

Consequently, the nonlinear PFs in Case I closely resemble the linear PFs, whereas the outcomes in Cases II and III exhibit more significant differences. Similar to Example 1, the responses of the $\lambda_1$ mode are reconstructed based on (20) in Figure 4. It's crucial to clarify that the scaling factors $\theta_i$ are distinct from the perturbation amplitudes $\alpha_i$, and these responses are utilized solely to illustrate the influence of $\theta_i$, which is similar to adjusting the perturbation amplitude $\alpha_i$ (as discussed in Remark 6).

This example indicates the significance of having unique PFs. In practical scenarios, mode shapes of a monitored nonlinear system can be obtained through signal processing techniques, such as Prony analysis (Xia et al., 2020), which measures system responses under small disturbances to approximate linear system behavior. However, obtaining a complete mode shape matrix $\mathbf{\Phi}$ can be challenging due to limitations in measuring devices. When calculating PFs, mode compositions cannot be directly derived from $\mathbf{\Phi}^{-1}$and must be obtained from the system model, where the choice of scaling factors for modes can have a significant impact.

Example 2 demonstrates that different choices of $\theta$-factors can lead to distinct values of nonlinear PFs even after normalization. Particularly, the third term in (21b) highlights that when the system exhibits noticeable nonlinearity, the impact of scaling becomes significant. According to the normal form theory, a near-resonance condition, even if it isn't a perfect resonance, has the potential to amplify nonlinearity (Dobson et al., 2001), thereby magnifying the influence of the scaling factor.

Certain literature, such as (Kundur, 1993; Sanchez-Gasca et al., 2005), recommends maintaining an inverse relationship between mode shapes and compositions (see Cases II and III in Example 2). Consequently, the selection of $\theta$-factors becomes critical to ensure $\boldsymbol{\psi}_i\boldsymbol{\phi}_i$ =1 for every mode $i$, while the values of $\xi$- or $\sigma$-factors can be relaxed.

### *4.4 Potential Future Work*

For a complex system, participation factors help identify the most important state variables for designing efficient control methods targeting dynamic modes of interest. When the system exhibits nonlinear dynamics, nonlinear participation factors become necessary, and establishing sufficient conditions for their uniqueness ensures the reliability of their values in identifying key variables for control. Consequently, future work will involve applying the sufficient conditions identified in this paper to control design. For instance, participation factors have demonstrated their effectiveness in model reduction and control of power systems with increasing penetration of inverter-based renewable energy resources (Xia, Ramasubramanian, et al, 2024; Sajjadi, 2022; Sajjadi, 2023; Sajjadi, 2024). When wide-area measurements, such as those from phasor measurement units, are available, participation factors may be estimated using data-driven algorithms (Xia, Yu, Sun, Shi & Huang, 2024). Therefore, the conclusions of this paper will aid in selecting and preprocessing measurement data for reliable estimation of participation factor estimation.

## 5. The uniqueness of other participation factors

Based on Theorems 1 and 2, this section proves the uniqueness conditions for other PFs, as detailed in Table 3.

**Corollary 1**: The PMISPF (Abed et al., 2000) are unique *if and only if the corresponding $\theta_i$ is unique*;

**Corollary 2**: The PSIMPF (Hashlamoun et al., 2009), Nonlinear PMISPF (Hamzi & Abed, 2020), and Modified PSIMPF (Iskakov, 2020) are unique *if all θ-factors are unique*.

**Table 3.** The variants of participation factors

| # | Notation | Full Name | Short name | Reference | Expression |
|---|---|---|---|---|---|
| 1 | $p_{ki}^{PMIS}$ | Probability Mode-in-state Participation Factor | PMISPF | (Abed et al., 2000) | (23) |
| 2 | $p_{ki}^{PSIM}$ | Probability State-in-mode Participation Factor | PSIMPF | (Hashlamoun et al., 2009) | (24) |
| 3 | $p_{ki}^{NPMIS}$ | Nonlinear Probability Mode-in-state Participation Factor | Nonlinear PMISPF | (Hamzi & Abed, 2020) | (25) |
| 4 | $p_{ki}^{MPSIN}$ | Modified Probability State-in-mode Participation Factor | Modified PSIMPF | (Iskakov, 2020) | (26) |
| 5 | $p_{ki}^{Data}$ | Data-driven Participation Factor | Data-driven PF | (Netto et al., 2019) | (27) |

*Proof*:

The traditional linear PF is independent of the selection of initial values. The PMISPF considers the quantity of the initial condition by computing a mode's average contribution to a state (Abed et al., 2000). Following a similar proof structure as in Theorem 1, (13) is applied to replace the $\boldsymbol{\phi}_i$ and $\boldsymbol{\psi}_i$ with $\hat{\boldsymbol{\phi}}_i$ and $\hat{\boldsymbol{\psi}}_i$, so we have

$$p_{ki}^{PMIS} = E\left\{\frac{\left(\boldsymbol{\psi}_i^{\mathrm{T}}\mathbf{x}_0\right)\phi_{ki}}{x_{k0}}\right\} = E\left\{\frac{\theta_i\left(\hat{\boldsymbol{\psi}}_i^{\mathrm{T}}\mathbf{x}_0\right)\hat{\phi}_{ki}}{\left(\cos\delta_i\right)x_{k0}}\right\}. \tag{23}$$

where $E\{\bullet\}$ represents the expectation operator, and $x_{k0}$ denotes the initial values for the $k$-th state variable in $\mathbf{x}$-space. It is obvious that $p_{ki}^{PMIS}$ in (23) is unique if and only if the corresponding $\theta_i$ is unique.

While distinguishing between state-in-mode and mode-in-state is unnecessary for a linear PF due to their identical nature, a study by (Hashlamoun et al., 2009)highlights that PSIMPF and PMISPF are not interchangeable. In this paper, another difference between them is exposed in the view of uniqueness:

$$
\begin{aligned}
p_{ki}^{PSIM} &= \begin{cases} E\left\{\dfrac{\psi_{ik}x_{k0}}{z_{i0}}\right\} & \lambda_i \in \mathbb{R} \\ E\left\{\dfrac{(\psi_{ik}+\psi_{ik}^{*})x_{k0}}{z_{i0}+z_{i0}^{*}}\right\} & \lambda_i \notin \mathbb{R} \end{cases} \\
&= \begin{cases} E\left\{\dfrac{\xi_i\hat{\psi}_{ik}x_{k0}}{\xi_i J_i(\theta_1,\theta_2,....\theta_n)}\right\} = E\left\{\dfrac{\hat{\psi}_{ik}x_{k0}}{J_i(\theta_1,\theta_2,....\theta_n)}\right\} & \lambda_i \in \mathbb{R} \\ E\left\{\dfrac{(\hat{\psi}_{ik}+\hat{\psi}_{ik}^{*})x_{k0}}{J_i(\theta_1,\theta_2,....\theta_n)+J_i^{*}(\theta_1,\theta_2,....\theta_n)}\right\} & \lambda_i \notin \mathbb{R} \end{cases}
\end{aligned} \tag{24}
$$

where $z_{i0} = \xi_i J_i(\theta_1, \theta_2, ..., \theta_n)$ denotes the initial values for the $i$-th state variable in **z**-space, $J_i$: $\mathbb{R}^n \rightarrow \mathbb{R}$ is a function of all $\theta$-factors, and * in the subscript shows its conjugate value (not conjugate function). Thus, the uniqueness of PSIMPF $p_{ki}^{PSIM}$ is related to all $\theta$-factors rather than $\theta_i$ in PMISPF.

The nonlinear PMISPF extends the concept of probability MISPF introduced in (23) by incorporating 2nd nonlinearity through normal form theory (Hamzi & Abed, 2020). For simplification, $\alpha_k$ = 1 is assumed in the following proof related to nonlinear PF. Based on Theorem 2, it is intuitive that it is unique if all $\theta$-factors are determined:

$$
p_{ki}^{NPMIS} = E\left\{\left.\frac{z_{i0}\phi_{ki}e^{\lambda_i t}}{x_k(t)}\right|_{t=0}\right\} = E\left\{\left.\frac{z_{i0}\left(\sigma_i\hat{\phi}_{ki}\right)e^{\lambda_i t}}{x_k(t)}\right|_{t=0}\right\} = E\left\{\left.\frac{\theta_i J_i(\theta_1,\theta_2,....\theta_n)\hat{\phi}_{ki}e^{\lambda_i t}}{\left(\cos\delta_i\right)x_k(t)}\right|_{t=0}\right\}. \tag{25}
$$

The modified PSIMPF extends its consideration to include the energy of the mode (Iskakov, 2020). Similar to (24), it is shown in (26) that all the $\theta$-factors appear in $p_{ki}^{MPSIM}$:

$$
\begin{aligned}
p_{ki}^{MPSIM} &= \frac{E\left\{(\psi_{ik}x_{k0})^{*}z_{i0}+z_{i0}^{*}(\psi_{ik}x_{k0})\right\}}{2E\left\{z_{i0}z_{i0}^{*}\right\}} \\
&= \frac{E\left\{\xi_i^2(\hat{\psi}_{ik}x_{k0})^{*}J_i(\theta_1,\theta_2,....\theta_n)+\xi_i^2 J_i^{*}(\theta_1,\theta_2,....\theta_n)(\hat{\psi}_{ik}x_{k0})\right\}}{2E\left\{\xi_i^2 z_{i0}(\theta_1,\theta_2,....\theta_n)J_i^{*}(\theta_1,\theta_2,....\theta_n)\right\}} \\
&= \frac{E\left\{(\hat{\psi}_{ik}x_{k0})^{*}J_i(\theta_1,\theta_2,....\theta_n)+J_i^{*}(\theta_1,\theta_2,....\theta_n)(\hat{\psi}_{ik}x_{k0})\right\}}{2E\left\{J_i(\theta_1,\theta_2,....\theta_n)J_i^{*}(\theta_1,\theta_2,....\theta_n)\right\}}
\end{aligned} \tag{26}
$$

□

In the data-driven PF approach described in (Netto et al., 2019), the PFs are determined using Koopman mode decomposition. Notably, Koopman modes are defined from a signal perspective, which may result in their mode shapes and composition coinciding with those defined in the system model. The definition of Koopman modes bears a resemblance to the structure of the PMISPF outlined in (23), unless a Koopman mode is specifically under consideration:

$$p_{ki}^{Data} = E\left\{\frac{\left(\mathbf{u}_i^{\mathrm{T}}\boldsymbol{\gamma}_0\right)v_{ki}}{\gamma_{k0}}\right\}, \tag{27}$$

where $\boldsymbol{\gamma}_0 \in \mathbb{R}^{n\times l}$ and $\gamma_{k0}$ are the initial values for the observable vector and $k$-th element. $\mathbf{u}_\mathrm{i}$ and $\mathbf{v}_\mathrm{i}$ are the $i$-th left and right eigenvector of the Koopman operator with

$$\begin{bmatrix}\mathbf{v}_1 & \ldots & \mathbf{v}_i & \ldots\end{bmatrix} = \mathbf{B}\begin{bmatrix}\mathbf{u}_1^{\mathrm{T}} & \ldots & \mathbf{u}_i^{\mathrm{T}} & \ldots\end{bmatrix}^{-1},$$

where $\mathbf{B} \in \mathbb{R}^{n\times l}$ is the matrix determined by the observed function and state variable. Note that (27) is similar to (23). Hence, if the scaling factors $\sigma_i^K$, $\xi_i^K$ and $\theta_i^K$ are introduced for Koopman mode by replacing $\boldsymbol{\phi}_i$ and $\boldsymbol{\psi}_i$ with $\mathbf{u}_\mathrm{i}$ and $\mathbf{v}_\mathrm{i}$ in (13). From (23) and (27), it is easily known that $p_{ki}^{Data}$ is unique if and only if the corresponding $\theta_i^K$ for Koopman mode $i$ is unique.

## 6. Conclusion

This paper has studied the conditions for uniqueness in various forms of PFs while considering scaling uncertainties in mode shapes and compositions. Three scaling factors were introduced, and the uniqueness of nonlinear PFs for linear or combination modes was thoroughly discussed, providing a sufficient condition. In contrast to perturbation amplitudes that impact nonlinear PFs from the perspective of state variables, scaling factors offer a means to adjust nonlinearity from the viewpoint of modes. Additionally, uniqueness conditions were established for several other PF variants. Understanding these sufficient conditions is crucial for applying the concept of PFs correctly in practical scenarios of stability analysis and control.

## Reference

Abed, E. H., Lindsay, D., & Hashlamoun, W. A. (2000). On participation factors for linear systems. *Automatica*, *36*(10), 1489-1496. https://doi.org/https://doi.org/10.1016/S0005-1098(00)00082-0

Amano, H., Kumano, T., & Inoue, T. (2006). Nonlinear stability indexes of power swing oscillation using normal form analysis. *IEEE Transactions on Power Systems*, *21*(2), 825-834. https://doi.org/10.1109/TPWRS.2006.873110

Arnold, V. I. (1988). *Geometrical Methods in the Theory of Ordinary Differential Equations* (J. Szücs, Trans.). Springer-Verlag New York.

Courtois, N., Klimov, A., Patarin, J., & Shamir, A. (2000, 2000//). Efficient Algorithms for Solving Overdefined Systems of Multivariate Polynomial Equations. (Ed.),^(Eds.). Advances in Cryptology — EUROCRYPT 2000, Berlin, Heidelberg.

Dobson, I., & Barocio, E. (2004). Scaling of normal form analysis coefficients under coordinate change. *IEEE Transactions on Power Systems*, *19*(3), 1438-1444. https://doi.org/10.1109/TPWRS.2004.831691

Dobson, I., & Barocio, E. (2005). Perturbations of weakly resonant power system electromechanical modes. *IEEE Transactions on Power Systems*, *20*(1), 330-337. https://doi.org/10.1109/TPWRS.2004.841242

Dobson, I., Zhang, J., Greene, S., Engdahl, H., & Sauer, P. W. (2001). Is strong modal resonance a precursor to power system oscillations? *IEEE Transactions on Circuits and Systems I: Fundamental Theory and Applications*, *48*(3), 340-349. https://doi.org/10.1109/81.915389

Dreesen, P., Batselier, K., & De Moor, B. (2012). Back to the Roots: Polynomial System Solving, Linear Algebra, Systems Theory. *IFAC Proceedings Volumes*, *45*(16), 1203-1208. https://doi.org/https://doi.org/10.3182/20120711-3-BE-2027.00217

Benallou, A., Mellichamp, D., & Seborg, D. (1983). On the number of solutions of multivariable polynomial systems. *IEEE Transactions on Automatic Control*, *28*(2), 224-227. https://doi.org/10.1109/TAC.1983.1103214

Garofalo, F., Iannelli, L., & Vasca, F. (2002). PARTICIPATION FACTORS AND THEIR CONNECTIONS TO RESIDUES AND RELATIVE GAIN ARRAY. *IFAC Proceedings Volumes*, *35*(1), 125-130. https://doi.org/https://doi.org/10.3182/20020721-6-ES-1901.00182

Hamzi, B., & Abed, E. H. (2020). Local modal participation analysis of nonlinear systems using Poincaré linearization. *Nonlinear Dynamics*, *99*(1), 803-811. https://doi.org/10.1007/s11071-019-05363-1

Hashlamoun, W. A., Hassouneh, M. A., & Abed, E. H. (2009). New Results on Modal Participation Factors: Revealing a Previously Unknown Dichotomy. *IEEE Transactions on Automatic Control*, *54*(7), 1439-1449. https://doi.org/10.1109/TAC.2009.2019796

Huang, Q., Wang, Z., & Zhang, C. (2009, 26-30 July 2009). Evaluation of the effect of modal interaction higher than 2nd order in small-signal analysis. (Ed.),^(Eds.). 2009 IEEE Power & Energy Society General Meeting.

Iskakov, A. B. (2020). Definition of State-In-Mode Participation Factors for Modal Analysis of Linear Systems. *IEEE Transactions on Automatic Control*, 1-1. https://doi.org/10.1109/TAC.2020.3043312

Kundur, P. (1993). *Power System Stability and Control*.

Liu, S., Messina, A. R., & Vittal, V. (2006). A Normal Form Analysis Approach to Siting Power System Stabilizers (PSSs) and Assessing Power System Nonlinear Behavior. *IEEE Transactions on Power Systems*, *21*(4), 1755-1762. https://doi.org/10.1109/TPWRS.2006.882456

Netto, M., Susuki, Y., & Mili, L. (2019). Data-Driven Participation Factors for Nonlinear Systems Based on Koopman Mode Decomposition. *IEEE Control Systems Letters*, *3*(1), 198-203. https://doi.org/10.1109/lcsys.2018.2871887

Perez-arriaga, I. J., Verghese, G. C., & Schweppe, F. C. (1982). Selective Modal Analysis with Applications to Electric Power Systems, PART I: Heuristic Introduction. *IEEE Transactions on Power Apparatus and Systems*, *PAS-101*(9), 3117-3125. https://doi.org/10.1109/TPAS.1982.317524

Samovol, V. S. (2004). Normal Form of Autonomous Systems with One Zero Eigenvalue. *Mathematical Notes*, *75*(5), 660-668. https://doi.org/10.1023/B:MATN.0000030974.08984.21

Sanchez-Gasca, J. J., Vittal, V., Gibbard, M. J., Messina, A. R., Vowles, D. J., Liu, S., & Annakkage, U. D. (2005). Inclusion of higher order terms for small-signal (modal) analysis: committee report-task force on assessing the need to include higher order terms for small-signal (modal) analysis. *IEEE Transactions on Power Systems*, *20*(4), 1886-1904. https://doi.org/10.1109/TPWRS.2005.858029

Shu, L., Messina, A. R., & Vittal, V. (2005). Assessing placement of controllers and nonlinear behavior using normal form analysis. *IEEE Transactions on Power Systems*, *20*(3), 1486-1495. https://doi.org/10.1109/TPWRS.2005.852052

Yang, D., & Sun, Y. (2022). SISO impedance-based stability analysis for system-level small-signal stability assessment of large-scale power electronics-dominated power systems. *IEEE Transactions on Sustainable Energy*, 13(1), 537–550. https://doi.org/10.1109/TSTE.2021.3119207

Yang, Z., Zhan, M., & Tang, W. (2023). Dominant device identification in multi-VSC systems by frequency-domain participation factor. *IEEE Access*, 11, 90190–90200. https://doi.org/10.1109/ACCESS.2023.3288148

Zhu, Y., Gu, Y., Li, Y., & Green, T. C. (2022). Participation analysis in impedance models: The grey-box approach for power system stability. *IEEE Transactions on Power Systems*, 37(1), 343–353. https://doi.org/10.1109/TPWRS.2021.3088345

Xue, Y., Feng, N., & Wang, Y. (2023). Resonance participation factor evaluation for renewable energy systems based on impedance scanning. 2023 IEEE 7th Conference on Energy Internet and Energy System Integration (EI2), 1858–1863. https://doi.org/10.1109/EI259745.2023.10512826

Zhan, Y., Xie, X., Liu, H., Liu, H., & Li, Y. (2019). Frequency-domain modal analysis of the oscillatory stability of power systems with high-penetration renewables. *IEEE Transactions on Sustainable Energy*, 10(3), 1534–1543. https://doi.org/10.1109/TSTE.2019.2900348

Songzhe, Z., Vittal, V., & Kliemann, W. (2001, 15-19 July 2001). Analyzing dynamic performance of power systems over parameter space using normal forms of vector fields. II. Comparison of the system structure. (Ed.),^(Eds.). 2001 Power Engineering Society Summer Meeting. Conference Proceedings (Cat. No.01CH37262).

Starrett, S. K., & Fouad, A. A. (1998). Nonlinear measures of mode-machine participation [transmission system stability]. *IEEE Transactions on Power Systems*, *13*(2), 389-394. https://doi.org/10.1109/59.667357

Tian, T., Kestelyn, X., Thomas, O., Amano, H., & Messina, A. R. (2018). An Accurate Third-Order Normal Form Approximation for Power System Nonlinear Analysis. *IEEE Transactions on Power Systems*, *33*(2), 2128-2139. https://doi.org/10.1109/TPWRS.2017.2737462

Tzounas, G., Dassios, I., & Milano, F. (2020). Modal Participation Factors of Algebraic Variables. *IEEE Transactions on Power Systems*, *35*(1), 742-750. https://doi.org/10.1109/TPWRS.2019.2931965

Wang, Z., & Huang, Q. (2017). A Closed Normal Form Solution Under Near-Resonant Modal Interaction in Power Systems. *IEEE Transactions on Power Systems*, *32*(6), 4570-4578. https://doi.org/10.1109/TPWRS.2017.2679121

Xia, T., Huang, & K., Sun, K. (2024). A Tensor Contraction-Based Approach for Efficient Computation of Nonlinear Participation Factors, *IEEE PES General Meeting, Seattle, WA, July 21-25* https://doi.org/10.1109/PESGM51994.2024.10688792

Xia, T., Yu, Z., Sun, K., Shi, D., & Huang, K. (2025). Estimation of participation factors for power system oscillation from measurements. IEEE Transactions on Industry Applications. https://doi.org/10.1109/TIA.2025.3530869

Wiggins, S. (2003). *Introduction to applied nonlinear dynamical systems and chaos*. Springer.

Xia, T., Yu, Z., Sun, K., Shi, D., & Wang, Z. (2020, 2-6 Aug. 2020). Extended Prony Analysis on Power System Oscillation Under a Near-Resonance Condition. 2020 IEEE Power & Energy Society General Meeting (PESGM).

Xia, T., & Sun, K. (2022). Time-variant nonlinear participation factors considering resonances in power systems. 2022 IEEE Power & Energy Society General Meeting (PESGM), 1–5. https://doi.org/10.1109/PESGM48719.2022.9916932

Xia, T., Ramasubramanian, D., Wang, W., Sun, K., & Huang, K. (2024). Grid-forming inverter placement for power systems with high inverter-based resource penetration based on participation factors. IEEE Energy Conversion Congress and Exposition (ECCE), Phoenix, AZ, October 20–24, 2024.

Sajjadi, M., Huang, K., & Sun, K. (2022). Participation factor-based adaptive model reduction for fast power system simulation. 2022 IEEE Power & Energy Society General Meeting (PESGM), Denver, CO, USA, 1–5. https://doi.org/10.1109/PESGM48719.2022.9916866

Sajjadi, M., Xia, T., Xiong, M., Sun, K., Hoke, A., Tan, J., Wang, B., (2023) Estimation of Participation Factors Using the Synchrosqueezed Wavelet Transform, *IEEE PES General Meeting*, Orlando, FL, July 16-20, 2023 ttps://doi.org/10.1109/PESGM52003.2023.10252202

Sajjadi, M., Xia, T., Xiong, M., Sun, K., Hoke, A., Wang, B., Tan, J. (2024), A Participation Factor-Based Approach for Defining the EMT Model Boundary for Power System Simulations with Inverter-Based Resources, *Annual Conference of the IEEE Industrial Electronics Society*, Chicago, IL, November 3 - 6, 2024.